# A numerical comment on the "tiny" oscillations and a heuristic "conjecture"

Danilo Merlini[1], Massimo Sala[2] and Nicoletta Sala[3]

[1] *CERFIM/ISSI, Locarno, Switzerland, e-mail: merlini@cerfim.ch*
[2] *Independent Researcher*
[3] *Institute for the Complexity Studies, Rome, Italy*

**Abstract**

In the first part of this short work (in the form of a comment) we add the plots of two more values of the Li-Keiper coefficients $\lambda_5$ and $\lambda_6$, computed as in our recent work where the first four values were in particular given.
This for the trend as well as for the oscillating path ("tiny", the term coined by Maslanka in his pioneering work). Then, in the second part looking at the tiny oscillations, we propose a "numerical conjecture" in a more strong form, i.e. with a logarithmic behaviour and carry out a short numerical experiment on the new "numerical conjecture".

**0. Short Introduction**

In a previous work [1], connected with other pioneering works on the subject [2, 3, 4, 5, 6, 7, 8] where the method of Baez-Duarte and the "strategy" of Malanska were considered, we have given the values of the first Li-Keiper coefficients (n up to n=4) involving values at integer or half-integer arguments of the Zeta function.
Here, as a comment, we first give the plots of two more calculated Li-Keiper coefficients. Then we look - considering recent advances in this direction – at a stronger – but more gratuitous numerical conjecture then the one we have advanced and discussed in [1] – and carry out a numerical experiment - leaving more elaborated experiments on our weaker conjecture [1] for two subsequent works on the subject [10, 11].

**1. Fluctuations: $\log((s-1)\cdot\zeta(s))$ with Pochammer's Polynomials, coefficients of $z^5$ and $z^6$ (z=1-1/s).**

We recall here that the log of the expansion of the Riemann $\xi$ Function is given by:

$$\log\left(\xi\cdot\left(\frac{1}{1-z}\right)\right) = \log\left[\frac{1}{2}\cdot\left(\frac{z}{(1-z)^2}\right)\cdot\pi^{-(1/(2\cdot(1-z)))}\cdot\Gamma(1/(2\cdot(1-z)))\cdot\zeta(1/(1-z))\right] \quad (1)$$

and

$$\log\left(\xi\cdot\left(\frac{1}{1-z}\right)\right) = \log\frac{1}{2} + \sum_{n=1}^{\infty}\lambda_n\cdot\frac{z^n}{n}$$

In Figures 1 and 2, the plots of the two coefficients (tiny fluctuations), coefficients of $\frac{z^5}{5}$ and of $\frac{z^6}{6}$ in the expansion of the tiny part $\log\left(\left(\frac{z}{1-z}\right)\zeta\cdot\left(\frac{1}{1-z}\right)\right)$.

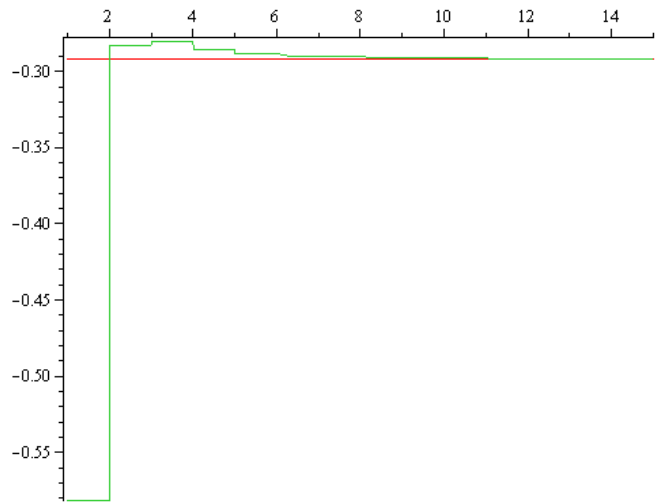

Fig. 1 The coefficient of $z^5$ in the expansion of $\log((s-1).\text{Zeta}(s))$
where $s = 1/(1-z)$, around $z = 0$ ($s=1$), in the range $1 < p < 15$.
(in red $-\lambda_t/5 = -(1.45826850020)/5 = -0.2916537000$
taken from Ref. [2]. Notice that $f(15) = -0.291599044$.
and that $f(1) = -0.58158...$

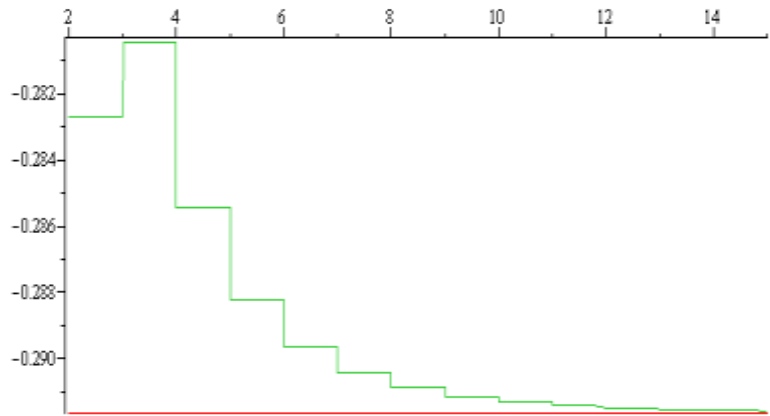

Fig. 2 The coefficient $-f(p)$, of $z^5$ in the expansion of $\log((s-1).\text{Zeta}(s))$
where $s = 1/(1-z)$, around $z=0$ ($s=1$), in the range $2 < p < 15$.
(in reed $-\lambda_t/5 = -(1.45826850020)/5 = -0.2916537000$
taken from Maslanka [2]. Notice that $f(23) = -0.291653084..$

**Remark**

Eventually our series in z so obtained [1] may be asymptotic and valid only up to some finite values of n in the coefficient of $z^n$. We have now additionally computed the coefficient of $z^5$ by means of our series with the Pochammer's Polynomials and found (See the above pictures) a value close to that of [Ref 2]. It should be noted that the function is not increasing in all domain but after p=4, decreasing monotonically to the expected value – 0.291.

Moreover the coefficient of $z^6$, that is $\lambda_t(n)$, (where t means tiny) has been also calculated for n=6 too and the result is given on the two Figures below.

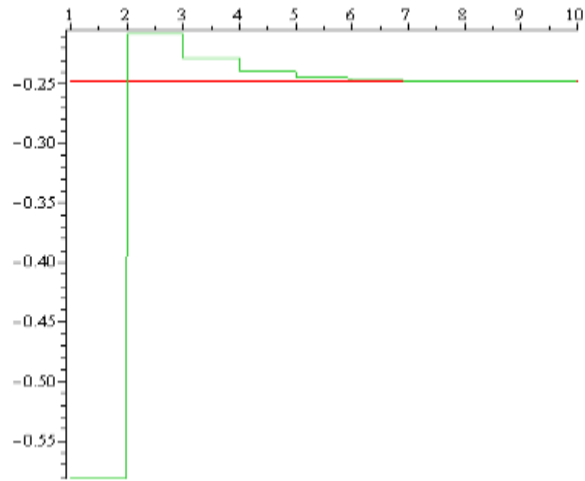

Fig. 3 The coefficient -f(p) of $z^6$ in the expansion of $\log((s-1)\cdot \text{Zeta}(s))$ (in red the value of Ref [2] , - 0.2480497212).

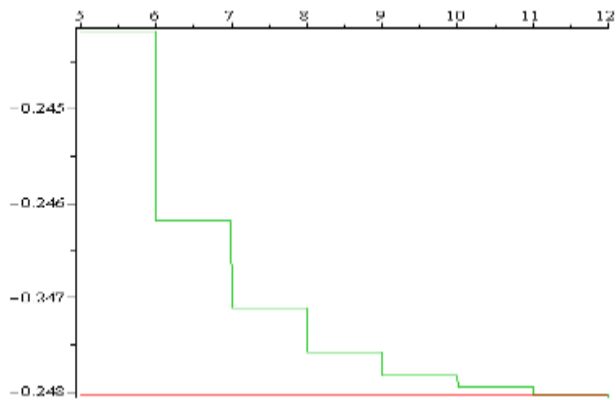

Fig. 4 The coefficient -f(p) of $z^6$ as above in the range 5<p<12 (in red the value of Ref [2] , Table in penultimate column of Maslanka [2] ($\lambda_{tiny}(6) = 1.48829832721 = 6\cdot 0.2480497212$).

**Comment**
The conjecture in our previous work on the subject, i.e. that the tiny fluctuations connected with the expansion of $\log((s-1)\cdot\zeta(s))$, i.e. of $\log([z/(1-z)]\cdot\zeta(s))$ is given by $|\lambda_{tiny(n)}| < f(1)\cdot n = 0.58158\cdot n$ may be strengthened heuristically, for example as $|\lambda_{tiny(n)}| < \gamma\cdot n$ where $\gamma = 0.577215..$ is the Euler-Mascheroni constant, taking into account the fact that the first 30 values of $\lambda_{tiny(n)}(n)/n$ are decreasing with n, as it may be seen from the Table below with a few of decimals. On the other hand, analyzing the numerical values at our disposal as reported below, we may strength the conjecture to $|\lambda_{tiny(n)}| = o(n) = a\cdot \log(n)$ with a~2.

| n | $\lambda_{tiny}(n)/n\cdot\gamma$ | $\lambda_{tiny}(n)$ | $2\cdot\log(n)$ |
|---|---|---|---|
| 1 | 1 | 0.577215 | |
| 2 | 0.837542 | 0.966885 | 1.386294 |
| 3 | 0.704934 | 1.220696 | 2.197 224 |
| 4 | 0.595786 | 1.375588 | 2.772588 |
| 5 | 0.505276 | 1.458268 | 3.218875 |
| 6 | 0.429734 | 1.488298 | 3.583518 |
| 7 | 0.366337 | 1.480190 | 3.891820 |
| 8 | 0.312893 | 1.444855 | 4.158883 |
| 9 | 0.267682 | 1.399596 | 4.394449 |
| 10 | 0.229342 | 1.323802 | 4.605170 |
| 15 | 0.108860 | 0.942358 | 5.416100 |
| 20 | 0.058093 | 0.670652 | 5.991464 |
| 25 | 0.052749 | 0.594962 | 6.437751 |
| 26 | 0.040167 | 0.602799 | 6.516193 |
| 27 | 0.039627 | 0.617452 | 6.591673 |
| 28 | 0.039511 | 0.638020 | 6.664409 |
| 29 | 0.039724 | 0.665174 | 6.734591 |
| 30 | 0.040235 | 0.697102 | 6.802394 |
| 31 | 0.041021 | 0.733544 | 6.867974 |

Maslanka gives also:

| n | $\lambda_{tiny}(n)/n\cdot\gamma$ | $\lambda_{tiny}(n)$ | $2\cdot\log(n)$ |
|---|---|---|---|
| 100 | 0.0108 | 0.628752 | 9.210 |
| 500 | 0.0092 | 2.663502 | 12.429 |
| 1000 | 0.0030 | 1.756264 | 13.815 |
| 2000 | 0.0093 | 10.76850 | 15.201 |
| 3000 | 0.0012 | -2.09000 | 16.012 |

From the plots in [5, 6] we read approximatively:

| n | $\lambda_{tiny}(n)/n\cdot\gamma$ | $\lambda_{tiny}(n)$ | $2\cdot\log(n)$ |
|---|---|---|---|
| 200 | 0.0311 | 3.600 | 10.596 |
| 760 | 0.0144 | -6.33 | 13.266 |
| 840 | 0.0173 | 8.4 | 13.466 |
| 900 | 0.0173 | 9.000 | 13.604 |
| 1870 | 0.0104 | -11.226 | 15.067 |
| 3300 | 0.00519 | -9.900 | 16.203 |
| 5080 (0.00335) | 0.00573 | 16.830 | 17.066 |
| 5500 | 0.00519 | 16.500 | 17.225 |
| 6500 | 0.00346 | 13.000 | 17.559 |
| 8000 | 0.00346 | 16.000 | 17.774 |

The above conjecture appears of course as gratuitous, but is as a preliminary preparation to a less heuristic treatment in a numerical as well as analytical context of our subsequent works in this direction concerning the tiny, the trend and the complete Li-Keiper coefficients: in those works we introduce and study a special approximation procedure [10, 11].

As a general comment it is important to add the following: assuming Riemann Hypothesis, Oesterlé [7] has shown that the fluctuations are o(n). Lagarias [8] has shown that o(n) (as fluctuation over the trend given by (n/2·log(n)+c·n), may be improved to be $o(n) = O(\sqrt{n} \cdot log n)$, while Reyna [6] has shown that o(n) = n·$y_n$ with the sequence {$y_n$} in $l_2$. In this comment, we have suggested - on the analysis of the known numerical results at our disposal - (See Table above) and Keiper [5] that the fluctuations may be conjectured to be even more small and given by a·log(n) with the number a around 2 for small n, but from recent results around n equal 80000, a should be at least 5 [12].
A Formula with three terms would be:

$$\lambda(n) = \lambda_{trend}(n) + \lambda_{tiny}(n) = (n/2) \cdot \log(n) - 1.13...(n) \pm a \cdot \log(n) \qquad (2)$$

A formula of this type in another range of z and in another variable (N) in the domain of absolute convergence and not related to the above one is given in the Appendix (variable N!).

**Conclusion**

In the first part of this comment we have computed two more coefficients (of $\frac{z^5}{5}$ and $\frac{z^6}{6}$) of the expansion at z=0 of the oscillating part and our previous conjecture may be "improved" to be γ·n, where γ is the Euler-Mascheroni constant. In the second part, we have carried out small numerical experiments on some known numerical results and we have advanced a more strong but more gratuitous conjecture that the tiny oscillations may be bounded in absolute value by a·log(n) (better $O(n^\varepsilon)$, $\forall \varepsilon > 0$) where n is the index n of the coefficient of $z^n$ for at large n.
To deepen some aspects of the research in course, see for example, Albeverio-Cebulla [13] and Albeverio-Cacciapuoti [14]. To deepen the general aspects of Riemann Hypothesis see Broughan [15, 16].

**Appendixes A and B.**

**A.** The three functions appearing in the proposed behavior are n·log(n), n, log(n), n·log(n)/log(n) = n and n/(n/log(n)) = log(n). Here, as a curiosity, we give explicitely a function of this type in a new variable denoted by N, (not to be confused with the index n of the λ's). We consider the function ξ in the domain of absolute convergence and ξ(s) for s→N.
Then

$$\log(\xi(N)) = \log((½) \cdot N \cdot (N-1) \cdot \pi^{(-N/2)} \cdot \Gamma(N/2) \cdot \varsigma(N)) =$$

$$= \log(N-1) \cdot \varsigma(N)) - (N/2) \cdot \log(\pi) + \log(\Gamma(1+N/2)) . \qquad (3)$$

$$\log[(N-1) \cdot \pi^{-(N/2)} \cdot \Gamma(1+N/2)] = -(N/2) \cdot \log(\pi) + \log(N) + \log(\Gamma(1+N/2))$$

for the last term, the Stirling approximation gives:

$$\log(\Gamma(1+N/2)) \sim (N/2) \cdot \log(N) + (N/2 \cdot (-\log(2)-1) + (1/2) \cdot \log(N)$$

Thus, for ξ(s), we have:

$$\log(\xi(N)) \sim (N/2) \cdot \log(N) + (N/2) \cdot (-\log(2\pi)-1) + (1/2) \cdot \log(N) + \log((N-1)..\varsigma(N))$$

$$\sim (N/2) \cdot \log(N) + (N/2) \cdot (-\log(2\pi)-1) + (3/2) \cdot \log(N) \qquad (4)$$

with the three terms "of some interest".

**B.** We now carry out a numerical experiment to check a sum on the Li-Keiper coefficients. Then, the Formula gives us:

$$\xi(1/(1-z)) = -\log(2) + \sum_{n=1}^{\infty}(1/n) \cdot \lambda_n \cdot z^n$$

and for z=1/2 it is equal to $\log(\pi/3) = 0.04611759699$ for z=1/2 (assuming on RH that the sum exists! here is to remember that $\varsigma(2) = \pi^2/6$, Euler). Then the left hand side of the above Equation, using the first 15 values of the Table (Ref. [2]) the contribution is: 0.04610606601.
From a n=16 to infinity we assume the formula $(n/2) \cdot \log(n) + c \cdot n + a \cdot \log(n)$, with $c = (1/2) \cdot (\gamma - \log(2\pi) - 1) = -1.13$...with **a** unknown, for $\lambda(n)$.

We find the approximate Equation:

0.04610606..+0.0007357866258 – 0.0005864142430+a·0.0008636699215- 0.046117597181290 =0.
With the solution a= -1.59599 = -1.596.
Then
$$\lambda(n) \sim (n/2) \cdot \log(n) + c \cdot n \pm 1.596 \cdot \log(n)..., \text{ with } c = (\gamma - \log(2\pi) - 1).$$

If -on the other hand -we would assume the behavior $a \cdot \sqrt{n} \cdot \log(n)$,
 (assuming the R.H. [8, 9]), we obtain the Equation:
0.04610606 +0.0007357866258-0.0005864142430+a·0.0003562045074-0.046117597181290=0,
a= -0.3869721386 = -0.386.

then: $\lambda(n) \sim (n/2) \cdot \log(n) + c \cdot n \pm 0.386 \cdot \sqrt{n} \cdot \log(n)$ ... for this possibility. (5)

A plot, as illustration for the reader is given below.

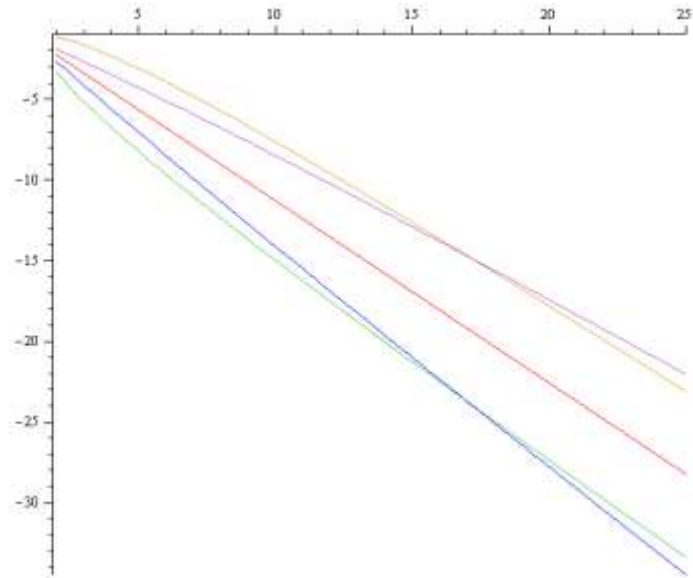

Fig. A.1 In red c·n ,in blue and violet c·n ±0.386 · √n ·log(n),in green and maroon  c·n ± 1.596·log(n) (around n=17, the plots intersect: we have not taken  the term (n/2)·log(n)).

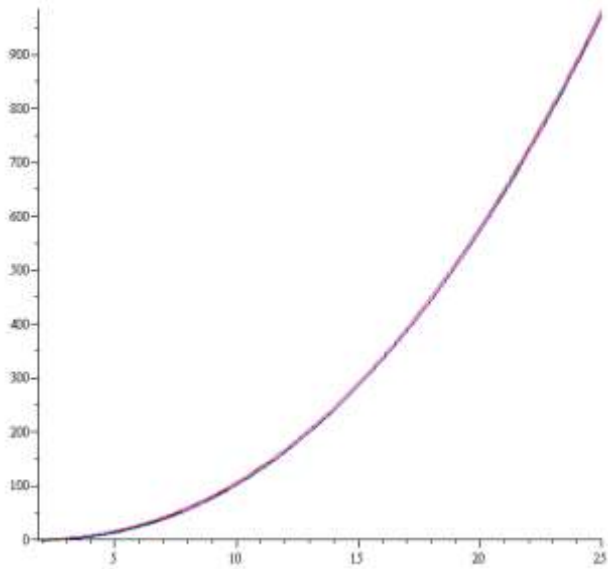

Fig. A.2 With the main term (n/2)·log(n), illustrating the tiny fluctuations.

## References

[1] D. Merlini, M. Sala, and N. Sala, Fluctuation around the Gamma function and a Conjecture, *IOSR International Journal of Mathematics*, 2019, volume 15 issue 1, 57-70

[2] K. Maslanka, "Baez-Duarte's criterion for the Riemann Hypothesis and Rice's Integrals", Math N.T. /0603713v2 /ArXiv (2018)

[3] L. Baez-Duarte: "A new necessary and sufficient condition for the Riemann Hypothesis": ArXiv: mat/0307215 (2003)

[4] Yu. V. Matiyasevich, Yet Another Representation for Reciprocals of the Nontrivial Zeros of the Riemann Zeta Function, *Mat. Zametki*, **97**:3 (2015), 471–474; *Math. Notes*, **97**:3 (2015), 476–479

[5] J.B. Keiper, Power Series Expansions of Riemann's $\xi$ Function. *Mathematics of Computation*, volume 58, number 198, April 1992, 765-773.

[6] J. Aria de Reyna, Asymptotics of Li-Keiper Coefficients (2011). Retrieved 10 January 2019: https://personal.us.es/arias/56-2RHB.pdf

[7] J. Oesterlé, J, Régions sans zéros de la fonction zêta de Riemann, typescript (2000, revised 2001, uncirculated).

[8] J.C. Lagarias, Li coefficients for automorphic L-functions. *Annales de l'Institut Fourier*, Volume 57 (2007) no. 5, pp. 1689-1740. doi: 10.5802/aif.2311. http://www.numdam.org/item/AIF_2007__57_5_1689_0/

[9] A. Voros, Sharpening of Li's Criterion for the Riemann Hypothesis, *Mathematical Physic, Analysis and Geometry* (2006) 9, issue 1, 53-63. https://doi.org/10.1007/s11040-005-9002-8

[10] D. Merlini, M. Sala, N. Sala, Some approximations on the Li-Keiper coefficients: numerical experiments, to be posted

[11] D. Merlini, M. Sala, N. Sala, Analysis of the Li-Keiper coefficients by a "perturbation" of the K function, in preparation.

[12] F. Johansson, Rigorous high-precision computation of the Hurwitz zeta function and its derivatives, arXiv:1309.2877v1 [cs.SC] (2013)

[13] S. Albeverio, C. Cebulla, Müntz formula and zero free regions for the Riemann zeta function, *Bull. Sci. math*. 131 (2007), 12–38

[14] S. Albeverio, C. Cacciapuoti, The Riemann zeta in terms of the dilogarithm, *Journal of Number Theory,* Volume 133, Issue 1 (2013) 242-277, https://doi.org/10.1016/j.jnt.2012.06.002

[15] K. Broughan, *Equivalents of the Riemann Hypothesis* (2017), vol. 1, Cambridge University Press.

[16] K. Broughan, *Equivalents of the Riemann Hypothesis* (2017), vol. 2, Cambridge University Press.